\documentclass[10pt]{amsart}
\usepackage{amssymb}

\def\today{\number\day\space\ifcase\month\or   January\or February\or
   March\or April\or May\or June\or   July\or August\or September\or
   October\or November\or December\fi\   \number\year}

\newcounter{TmpEnumi}

\theoremstyle{definition}
\newtheorem{thm}{Theorem}[section]
\newtheorem{lem}[thm]{Lemma}
\newtheorem{prp}[thm]{Proposition}
\newtheorem{dfn}[thm]{Definition}
\newtheorem{cor}[thm]{Corollary}

\newtheorem{rmk}[thm]{Remark}

\newcommand{\beq}{\begin{equation}}
\newcommand{\eeq}{\end{equation}}
\newcommand{\beqr}{\begin{eqnarray*}}
\newcommand{\eeqr}{\end{eqnarray*}}
\newcommand{\bal}{\begin{align*}}
\newcommand{\eal}{\end{align*}}
\newcommand{\bei}{\begin{itemize}}
\newcommand{\eei}{\end{itemize}}

\newcommand{\tie}{}

\newcommand{\af}{\alpha}
\newcommand{\bt}{\beta}
\newcommand{\gm}{\gamma}
\newcommand{\dt}{\delta}
\newcommand{\ep}{\varepsilon}

\newcommand{\io}{\iota}

\newcommand{\kp}{\kappa}
\newcommand{\ph}{\varphi}
\newcommand{\ps}{\psi}
\newcommand{\rh}{\rho}

\newcommand{\ta}{\tau}

\newcommand{\C}{{\mathbb{C}}}
\newcommand{\N}{{\mathbb{N}}}

\newcommand{\Asym}[1]{{C_{\mathrm{b}} ([1, \infty), \, #1)
                           / C_0 ([1, \infty), \, #1)}}
\newcommand{\Bcf}[1]{{C_{\mathrm{b}} ([1, \infty), \, #1)}}
\newcommand{\btu}{\bigtriangleup}
\newcommand{\conv}{{\mathrm{conv}}}
\newcommand{\id}{{\mathrm{id}}}

\newcommand{\spn}{{\mathrm{span}}}

\newcommand{\andeqn}{\,\,\,\,\,\, {\mbox{and}} \,\,\,\,\,\,}

\newcommand{\ca}{C*-algebra}
\newcommand{\ct}{continuous}
\newcommand{\pj}{projection}
\newcommand{\hm}{homomorphism}

\newcommand{\Wolog}{Without loss of generality}

\newcommand{\mops}{mutually orthogonal \pj s}

\newcommand{\cfn}{continuous function}

\title[The Calkin algebra has outer automorphisms]{The
       Calkin algebra has outer automorphisms}

\author{N.~Christopher Phillips}

\author{Nik Weaver}

\date{22~June 2006}

\address{Department of Mathematics, University of Oregon,
       Eugene OR 97403-1222, USA.}

\email[]{ncp@darkwing.uoregon.edu}

\address {Department of Mathematics,
          Washington University,
          Saint Louis MO 63130, USA}

\email[]{nweaver@math.wustl.edu}

\thanks{Research of the first author
  partially supported by NSF grant DMS~0302401.}

\begin{document}

\begin{abstract}
Assuming the continuum hypothesis,
we show that the Calkin algebra has $2^{\aleph_1}$
outer automorphisms.
\end{abstract}

\maketitle

\indent
Let $H$ be a separable infinite dimensional Hilbert space,
let $L (H)$ be the
algebra of bounded operators on $H$,
let $K (H)$ be the algebra of compact operators on $H$,
and let $Q = L (H)/K (H)$ be the Calkin algebra.
A long-standing problem asks whether every automorphism
of $Q$ is inner,
that is, of the form $x \mapsto u^* x u$ for some unitary $u \in Q$.

The earliest references we have found
to this problem
are implicitly in~\cite{BDF0} (see page~126)
and explicitly in~\cite{BDF} (see Remark~1.6(2)).
The automorphisms whose existence we prove are,
however, implemented by a unitary on every separable subalgebra
of the Calkin algebra,
and are thus not interesting from the point of view of extensions.
In particular, it remains open whether
there exists an automorphism of the Calkin algebra
which sends the image of the unilateral shift to its adjoint.

It is of course well known that, for any Hilbert space $H,$
all automorphisms of $L (H)$ are inner.
Corollary~8.8 of~\cite{IPP}
provides a factor of type~II$_1$ with separable predual
such that all automorphisms are inner.
(In this context, recall that,
by Corollary~5.13 of~\cite{SZ},
all C*~automorphisms of a von Neumann algebra
are in fact von Neumann algebra automorphisms.)

If $A$ is a separable simple \ca\  such that every automorphism of $A$
is inner (in the multiplier algebra if $A$ is not unital),
then $A$ must be isomorphic to the algebra of compact operators
on some Hilbert space.
This is Corollary~3 of~\cite{Sk},
but, as pointed out by George Elliott,
can also be derived from the existence of
nontrivial central sequences~\cite{AP2}.

Without simplicity, this is not true.
In fact, it fails even for commutative \ca s.
There exists a compact metric space $X,$
with more than one point,
such that the only homeomorphism $h \colon X \to X$ is the identity.
We are grateful to Greg Kuperberg for providing this example.
We don't know a reference,
but the construction of such a subset of the plane is easy to outline.
Start with a line segment in the plane.
Attach another line segment to the midpoint.
Now there are three line segments, the new one and two
halves of the original.
Attach two more line segments to one midpoint,
three more to another, and four more to the last,
ensuring that no new line segment meets the old set at any other point.
Order the midpoints of all the line segments now present,
attach five line segments to the first, six to the next, etc.
Repeat infinitely often,
and take the closure.
If each new line segment is chosen short enough,
the resulting space will be a compact contractible subset of the plane.
The set of midpoints, at all stages combined, is dense.
Removal of any one of them leaves a set with a finite
number of connected components,
different for each one, so any homeomorphism fixes every point
of this dense subset.

The Calkin algebra may be regarded as a noncommutative or ``quantum''
analog of the Stone-\v{C}ech remainder
$\N^* = \beta \N \setminus \N$ \cite{W},
and the question whether there are nontrivial automorphisms of
$Q$ is a natural analog of the question whether there are
nontrivial self-homeomorphisms of $\N^*$.
(A self-homeomorphism
of $\N^*$ is {\emph{trivial}} if it
arises from an ``almost permutation'' of $\N$,
that is, a bijection
between two cofinite subsets of $\N$.)
It is known that the
continuum hypothesis implies that nontrivial self-homeomorphisms
of $\N^*$ do exist~\cite{R}.
(However, it is consistent with ZFC that there are
no nontrivial self-homeomorphisms.
See Chapter~IV of~\cite{Sh}
for the original proof,
and see \cite{SS}
and~\cite{Vl} for more recent, shorter, proofs.)
We show here that the basic structure of the existence proof
can be adapted to establish the analogous result for~$Q$.

In C*-algebra language, the argument for $\N^*$ goes roughly
like this.
Assuming the continuum hypothesis, we can express $C(\N^*)$,
the algebra of continuous functions on $\N^*$,
as the union of a nested transfinite sequence of separable
unital C*-subalgebras $A_{\af}$ for $\alpha < \aleph_1$.
The algebra $C (\N^*)$ has the following injectivity property:
if $A$ is a separable unital
commutative C*-algebra and $B$ is a unital C*-subalgebra of $A$ then
every unital $*$-monomorphism from $B$ into $C (\N^*)$
extends to a $*$-monomorphism from $A$ into $C (\N^*)$.
Using this fact we are able to build a $*$-monomorphism from
$C (\N^*)$ into itself by defining it sequentially on the subalgebras
$A_{\af}$,
and using a back-and-forth construction we can ensure that
the resulting map is surjective, hence an automorphism.
Now for every
$\alpha$ there are at least two ways of extending a $*$-monomorphism
$A_{\af} \to C(\N^*)$ to $A_{\alpha + 1}$, so the construction
produces $2^{\aleph_1}$ distinct automorphisms of $C(\N^*)$,
or equivalently, $2^{\aleph_1}$ distinct self-homeomorphisms of $\N^*$.
But (again using the continuum hypothesis)
there are only $\aleph_1$ trivial self-homeomorphisms,
so some self-homeomorphisms must be nontrivial.

This proof technique does not work straightforwardly for the Calkin
algebra because $Q$ does not have the analogous injectivity property.
Let $s \in Q$ be the image of the unilateral shift in $L (H).$
Let $B = C (S^1),$ let $u$ be its standard unitary generator,
and let $\ph \colon B \to Q$ be the \hm\  such that $\ph (u) = s.$
Let $A = C ([0, 1]),$ and let $\io \colon B \to A$
be $\io (f) (t) = f ( \exp ( 2 \pi i t)).$
Then there is no \hm\  $\ps \colon A \to Q$ such that
$\ps \circ \io = \ph.$
For example, $\io (u)$ has a square root but $\ph (u) = s$
does not.
(This example is based on a suggestion of John McCarthy.)


Instead, we build up our automorphism of $Q$ along a nested
transfinite sequence of separable subalgebras $A_{\af}$ by
constructing a sequence of unitaries $u_{\af} \in Q$ which
implement automorphisms of the $A_{\af}$.
It is trivial to
extend such an automorphism of $A_{\af}$ to one of $A_{\alpha + 1}$
since we can just take $u_{\alpha + 1} = u_{\af}$, and using
Voiculescu's double commutant theorem
(Corollary~1.9 of~\cite{V}),
it is not hard to
see that this can always be done in more than one way.
The difficulty
appears at limit stages, where we have to find a unitary which
implements an automorphism of $A_{\af}$ that agrees with
the automorphisms of $A_{\bt},$ for $\beta < \alpha$,
that have already been defined.
We accomplish this using techniques developed
by Manuilov and Thomsen~\cite{MT}, extending previous work by
Voiculescu, which allow us to realize an asymptotically inner
automorphism as an inner automorphism.

The principal technical difficulty arises in ensuring that the
new unitaries which appear at limit stages are linked to previous
unitaries by paths of bounded length.
This is needed to keep the
transfinite induction going because we require asymptotically inner
automorphisms at future limit stages.
We develop the required tools in the first section.
The second section contains the proof of the main theorem.


We adopt the conventions throughout this paper
that subalgebras of \ca s are assumed closed and selfadjoint,
and that \hm s are $*$-\hm s.

We owe thanks to a number of people for helpful discussions, but
especially to Charles Akemann, Don Hadwin,
Klaus Thomsen, and Eric Wofsey.

\section{Homomorphisms to outer multiplier algebras}\label{Sec:ToOutM}

\indent
The main result of this section is Proposition~\ref{P:RelMV},
although it is Corollary~\ref{P:RelMVb} that we actually use later.
Proposition~\ref{P:RelMV} is a relative version
of the construction of Section~3 of~\cite{MT},
which, starting from an asymptotic morphism from
a separable \ca\  to an outer multiplier algebra,
produces a true \hm.
In this context, ``relative'' means that if the restriction
of the asymptotic morphism to a subalgebra comes from a \hm,
then the \hm\  we construct can be chosen to agree with that
\hm\  on the subalgebra.

We thank Klaus Thomsen for providing essential help with the material
in this section.

Recall the Bartle-Graves Selection Theorem~\cite{BG}:

\begin{thm}\label{T:BGS}
Let $E$ and $F$ be Banach spaces, and let
$q \colon E \to F$ be a \ct\  surjective linear map.
Then there exists a \cfn\  $g \colon F \to E$
(not necessarily linear)
such that $q \circ g = \id_{F}.$
\end{thm}

\begin{proof}
Let $X$ be the vector space of all \cfn s from $F$ to $E,$
let $Y$ be the vector space of all \cfn s from $F$ to $F,$
and define ${\overline{q}} \colon X \to Y$
by ${\overline{q}} (f) = q \circ f$.
Theorem~4 of~\cite{BG} implies that ${\overline{q}}$ is
surjective.
Therefore we can find $g \in X$
such that ${\overline{q}}(g) = \id_{F},$ as desired.
\end{proof}

We need a relative version of this theorem.

\begin{lem}\label{L:RelBGS}
Let $E$ and $F$ be Banach spaces, and let
$q \colon E \to F$ be a \ct\  surjective linear map.
Let $M \subset E$ and $N \subset F$ be closed subspaces
such that $q (M) = N,$
and let $f \colon N \to M$ be a \cfn\  such that
$q \circ f = \id_{N}.$
Then there exists a \cfn\  $g \colon F \to E$
such that $q \circ g = \id_{F}$
and $g |_N = f.$
\end{lem}

\begin{proof}
We first construct a retraction $h \colon F \to N,$
that is, a \cfn\  $h \colon F \to F$ (not necessarily linear)
such that $h (\xi) \in N$ for all $\xi \in F$
and $h (\xi) = \xi$ for all $\xi \in N.$
To get it,
apply Theorem~\ref{T:BGS}
to the surjective Banach space map $p \colon F \to F / N$
to obtain a \cfn\  $k \colon F / N \to F$
such that $p \circ k = \id_{F / N}.$
Define $h (\xi) = \xi - k (p (\xi)) + k (0).$
One checks immediately that $h (\xi) = \xi$ for $\xi \in N$
and that $p (h (\xi)) = 0$ for all $\xi \in F.$

Now apply Theorem~\ref{T:BGS}
to get a \ct\  section $g_0 \colon F \to E.$
Define $g \colon F \to E$ by
$g (\xi) = g_0 (\xi - h (\xi)) + f (h (\xi)) - g_0 (0).$
Clearly $g$ is \ct,
and
\[
(q \circ g) (\xi)
   = (q \circ g_0) (\xi - h (\xi))
     + (q \circ f \circ h) (\xi)
     - (q \circ g_0) (0)
   = (\xi - h (\xi)) + h (\xi) - 0
   = \xi.
\]
If $\xi \in N$ then $g(\xi) = g_0(\xi - \xi) + f(\xi) - g_0(0) = f(\xi)$.
This completes the proof.
\end{proof}

\begin{lem}\label{L:Extend}
Let $B$ and $D$ be \ca s.
Let $A$ be a subalgebra of $B,$ and
let $J$ be an ideal in $D.$
Let
\[
\pi \colon D \to D/J
\andeqn
\kp \colon \Bcf{D/J} \to \Asym{D/J}
\]
be the quotient maps.
Let
$\ph \colon B \to \Asym{D/J}$
be a \hm,
and suppose that $\ph |_A$
factors through a \hm\  $\rh \colon A \to D/J.$
Then there exists a family of functions $L = (L_t)_{t \in [1, \infty)}$
from $B$ to $D$ satisfying the following conditions:
\begin{enumerate}
\item\label{LE:Ct}
$(t, b) \mapsto L_t (b)$ is jointly \ct.
\item\label{LE:Equi}
The family $(L_t)_{t \in [0, \infty)}$ is equicontinuous.
\item\label{LE:Bdd}
$t \mapsto L_t (b)$ is bounded for every $b \in B.$
\item\label{LE:Lift}
For $b \in B,$ if $y \in \Bcf{D/J}$
is the function $t \mapsto \pi ( L_t (b)),$
then $\kp (y) = \ph (b).$
\item\label{LE:OnA}
For every $a \in A,$
the function $t \mapsto L_t (a)$ is constant.
\end{enumerate}
\end{lem}

The maps $L_t$ are not required to be linear.

If $J = 0,$
then $L$ is an equicontinuous asymptotic morphism from
$B$ to $D$ such that for each $a \in A,$
the function
$t \mapsto L_t (a)$ is the constant function $t \mapsto \rh (a).$

\begin{proof}[Proof of Lemma~\ref{L:Extend}]
We also write $\pi$ for the map
\[
\Bcf{D} \to \Bcf{D/J}.
\]
Let $M \subset \Bcf{D}$ be the subspace consisting of all constant
functions,
and let $N \subset \Asym{D/J}$ be the image of the constant functions
in $\Bcf{D/J}.$
Then $\kp \circ \pi (M) = N,$
so Theorem~\ref{T:BGS} provides a \ct\  section (not necessarily linear)
$N \to M.$
Since
\[
\kp \circ \pi \colon \Bcf{D} \to \Asym{D/J}
\]
is surjective, Lemma~\ref{L:RelBGS}
now provides a \cfn\  %
\[
S \colon \Asym{D/J} \to \Bcf{D}
\]
such that $\kp \circ \pi \circ S$ is the identity on $\Asym{D/J}$
and such that if $x \in \kp (\rh (A))$ then $S (x) \in \Bcf{D}$
is a constant function.
For $b \in B$ and $t \in [1, \infty),$ define
$L_t (b) = S (\ph (b)) (t).$

Conditions~(\ref{LE:Bdd}), (\ref{LE:Lift}),
and~(\ref{LE:OnA}) are immediate.
We prove Condition~(\ref{LE:Equi}).
So let $b \in B$ and let $\ep > 0.$
Choose $\dt > 0$ such that whenever $x \in D/J$ satisfies
$\| x - \ph (b) \| < \dt,$
then $\| S (x) - S ( \ph (b)) \| < \ep.$
If now $c \in B$ satisfies $\| c - b \| < \dt,$
then for all $t \in [1, \infty)$ we have
\[
\| L_t (c) - L_t (b) \|
  = \| S (\ph (c)) (t) - S ( \ph (b)) (t) \|
  \leq \| S (\ph (c)) - S ( \ph (b)) \|
  < \ep.
\]
Condition~(\ref{LE:Ct})
follows from equicontinuity of $(L_t)_{t \in [1, \infty)}$
and continuity
of $t \to L_t (b)$ for each $b \in B.$
\end{proof}

\begin{prp}\label{P:RelMV}
Let $B$ be a separable \ca, let $A \subset B$ be a subalgebra,
and let $E$ be a nonunital \ca\  with a countable approximate identity.
Let $\mu \colon B \to \Asym{M (E) / E}$
be a \hm,
and suppose that $\mu |_A$
factors through a \hm\  $\rh \colon A \to M (E) / E.$
Then there exists a \hm\  $\ps \colon B \to M (E) / E$
such that $\ps |_A = \rh.$
\end{prp}

The \hm\  $\ps$ will be a ``folding'' of an asymptotic
morphism from $B$ to $M (E) / E$ obtained from $\ph,$
as in Section~3 of~\cite{MT}.
See especially Lemma~3.5 and the discussion after Remark~3.6
in~\cite{MT}.
For our construction, we need a slight strengthening
of Lemma~3.2 of~\cite{MT}.
Let $\pi \colon M (E) \to M (E) / E$ be the quotient map.
See the beginning of Section~3 of~\cite{MT}
for the definition of a unit sequence,
and see the discussion before Lemma~3.2 of~\cite{MT}
for the definition of a compatible pair.

\begin{lem}\label{L:MT3.2}
Let $B$ be a separable \ca,
let $A \subset B$ be a subalgebra,
and let $E$ be a nonunital \ca\  with a countable approximate identity.
Let $(L_t)_{t \in [1, \infty)}$
be an equicontinuous family of functions from $B$ to $M (E),$
such that $(\pi \circ L_t)_{t \in [1, \infty)}$
is an asymptotic morphism from $B$ to $M (E) / E,$
and such that $t \mapsto L_t (a)$
is a constant function for every $a \in A.$
Let $(v_n)_{n = 0}^{\infty}$ be a unit sequence in $E$ with $v_0 = 0.$
Then there exists a unit sequence $(u_n)_{n = 0}^{\infty}$ in $E$
such that:
\begin{enumerate}
\item\label{LK:Cnv}
For every $n,$
$u_n$ is in the convex hull $\conv ( \{ v_j \colon j \geq n \} ).$
\item\label{LK:Cmp}
$\big( (L_t)_{t \in [1, \infty)}, \, (u_n)_{n = 0}^{\infty} \big)$
is a compatible pair for $(\pi \circ L_t)_{t \in [1, \infty)}.$
\item\label{LK:OnA}
There is a dense subset $S \subset A$ such that
\[
\sum_{n = 1}^{\infty}
   \big\| (u_n - u_{n - 1})^{1/2} L_1 (a)
            - L_1 (a) (u_n - u_{n - 1})^{1/2} \big\|
   < \infty
\]
for every $a \in S.$
\end{enumerate}
\end{lem}

\begin{proof}
Choose countable dense sets
$\{ b_1, b_2, \ldots \}$ in the unit ball of $B$
and $\{ a_1, a_2, \ldots \}$ in the unit ball of $A.$
Choose $\dt_n > 0$ such that
whenever $s, t \in [1, \, n + 3]$ satisfy $| s - t | < \dt_n$
and whenever $1 \leq k \leq n,$
then $\| L_s (b_k) - L_t (b_k) \| < 2^{-n}.$
Choose finite $\dt_n$-dense sets $S_n \subset [1, \, n + 3]$
such that $S_1 \subset S_2 \subset \cdots.$
Define a finite subset $F_n \subset M (E)$ by
\[
F_n = \{ L_1 (a_1), L_1 (a_2), \ldots, L_1 (a_n) \}
      \cup \big\{ L_{s} (b_k)
              \colon {\mbox{$s \in S_n$ and
                              $1 \leq k \leq n$}} \big\}.
\]
Choose $\ep_n > 0$ with $\ep_n < 2^{-n},$
and also
(using polynomial approximations to the functional calculus)
so small that
whenever $D$ is a \ca\  and $a, x \in D$ satisfy
$0 \leq a \leq 1,$ $\| x \| \leq 1,$ and $\| a x - x a \| < 2 \ep_n,$
then $\big\| a^{1/2} x - x a^{1/2} \big\| < 2^{-n}.$
\Wolog\  $\ep_1 \geq \ep_2 \geq \cdots.$

We now construct $u_n$ by induction,
such that~(\ref{LK:Cnv}) holds, such that $u_n u_{n - 1} = u_{n - 1},$
and such that $\| u_n x - x u_n \| < \ep_n$ for all $x \in F_n.$
Take $u_0 = 0.$
Given $u_n,$ by~(\ref{LK:Cnv}) there is $N$ such that
$u_n \in \spn ( \{ v_0, v_1, \ldots, v_N \} ).$
\Wolog\  $N \geq n.$
The set $\conv ( \{ v_j \colon j \geq N + 1 \} )$
is a convex approximate identity for $E.$
(In particular, it is directed:
if $e_1, e_2, \ldots, e_l \in \conv ( \{ v_j \colon j \geq N + 1 \} ),$
then there is $m$ such that
$e_1, e_2, \ldots, e_l
   \in \conv ( \{ v_j \colon N + 1 \leq j \leq m \} ),$
and $v_{m + 1} \geq e_1, e_2, \ldots, e_l.$)
Using the lemma on Page~330 of~\cite{Ar}
and a direct sum trick,
as done in the proof of Theorem~1 of~\cite{Ar},
there exists $u_{n + 1} \in \conv ( \{ v_j \colon j \geq N + 1 \} )$
such that $\| u_{n + 1} x - x u_{n + 1} \| < \ep_{n + 1}$
for all $x \in F_{n + 1}.$
Because $k \leq N < j$ implies $v_j v_k = v_k,$
it follows that $u_{n + 1} u_n = u_n.$
Similarly $u_{n + 1} v_n = v_n.$
This completes the induction.

Since $(v_n)_{n = 0}^{\infty}$ is an approximate identity
and $u_{n + 1} \geq v_n,$
it follows that $(u_n)_{n = 0}^{\infty}$ is an approximate identity.

For $n \geq k,$ we have
\[
\| u_n L_1 (a_k) - L_1 (a_k) u_n \| < \ep_n
\andeqn
\| u_{n + 1} L_1 (a_k) - L_1 (a_k) u_{n + 1} \|
  < \ep_{n + 1} \leq \ep_n,
\]
so
\[
\big\| (u_{n + 1} - u_n)^{1/2} L_1 (a_k)
             - L_1 (a_k) (u_{n + 1} - u_n)^{1/2} \big\|
  < 2^{-n}.
\]
This implies Condition~(\ref{LK:OnA}).

For Condition~(\ref{LK:Cmp}),
let $b \in B$ and let $\ep > 0.$
We find $N$ such that $n \geq N$ implies that
$\sup_{t \in [1, \, n + 3]} \| u_n L_t (b) - L_t (b) u_n \| < \ep.$
Using equicontinuity, choose $\dt > 0$ such that whenever $c \in B$
satisfies $\| c - b \| < \dt$ and $t \in [0, \infty),$
then $\| L_t (c) - L_t (b) \| < \frac{1}{5} \ep.$
Choose $k$ such that $\| b_k - b \| < \dt.$
Choose $N \geq k$ and also so large that
$\ep_n < \frac{1}{5} \ep$ and $2^{-n} < \frac{1}{5} \ep.$
Now let $t \in [1, \, n + 3],$
and choose $s \in S_n$ such that $| s - t | < \dt_n.$
Then
\begin{align*}
& \| u_n L_t (b) - L_t (b) u_n \|
   \\
& \mbox{} \hspace{2em}
    \leq 2 \| L_t (b) - L_t (b_k) \|
            + 2 \| L_t (b_k) - L_s (b_k) \|
            + \| u_n L_s (b_k) - L_s (b_k) u_n \|  \\
& \mbox{} \hspace{2em}
    < 2 \big( \tfrac{1}{5} \ep \big) + 2 \cdot 2^{-n} + \ep_n.
\end{align*}
So
\[
\sup_{t \in [1, \, n + 3]} \| u_n L_t (b) - L_t (b) u_n \|
 \leq \tfrac{2}{5} \ep + 2 \cdot 2^{-n} + \ep_n
 < \ep,
\]
as desired.

The definition of a unit sequence at the beginning of Section~3
of~\cite{MT} requires one final condition, which we have not
verified: that there exist a strictly positive element $x \in E$
such that every $u_n$ can be obtained from $x$ by suitable
functional calculus.
But this is automatic using just $u_{n + 1} u_n = u_n$
and the fact that $(u_n)_{n = 0}^{\infty}$ is an approximate identity:
take $x = \sum_{n = 0}^{\infty} 2^{-n} u_n.$
\end{proof}

\begin{proof}[Proof of Proposition~\ref{P:RelMV}]
Let $\pi \colon M (E) \to M (E) / E$ be the quotient map.
Let $(L_t)_{t \in [1, \infty)}$ be as in Lemma~\ref{L:Extend},
with $M (E)$ in place of $D,$
with $E$ in place of $J,$
and with $\mu$ in place of $\ph.$
Set $\ph_t = \pi \circ L_t,$
giving an equi\ct\  asymptotic morphism from $B$ to $M (E) / E.$
We follow Lemmas~3.2 through~3.5 of~\cite{MT},
and the associated discussion, except that we substitute
Lemma~\ref{L:MT3.2} for Lemma~3.2 of~\cite{MT},
and we use the equi\ct\  lift $(L_t)_{t \in [1, \infty)}$ of $\ph.$
Let $\ps = \ph^{\mathrm{f}}$ be the folding of $\ph$ obtained from
Lemma~3.5 of~\cite{MT};
taking $u_0 = 0$
and with a suitable sequence $(t_n)_{n = 1}^{\infty}$ in $[1, \infty),$
it is a \hm\  $B \to M (E) / E$ given by
$\psi (b) = \pi (T (b))$ with
\[
T (b)
  = \sum_{n = 0}^{\infty}
    (u_n - u_{n - 1})^{1/2} L_{t_n} (b) (u_n - u_{n - 1})^{1/2},
\]
with convergence in the strict topology of $M (E).$

We need only prove that
$\psi (a) = \rh (a)$ for $a \in A.$
It suffices to prove this for $a$ in a dense subset of $A,$
and we use the subset of Lemma~\ref{L:MT3.2}(\ref{LK:OnA}).
Again with convergence in the strict topology,
we have
$L_1 (a) = \sum_{n = 0}^{\infty} (u_n - u_{n - 1}) L_1 (a),$
so
\[
T (a) - L_1 (a)
  = \sum_{n = 0}^{\infty} (u_n - u_{n - 1})^{1/2}
      \big[ L_1 (a) (u_n - u_{n - 1})^{1/2}
           - (u_n - u_{n - 1})^{1/2} L_1 (a) \big].
\]
For $a \in S,$
this series converges in norm
by Condition~(\ref{LK:OnA}) in Lemma~\ref{L:MT3.2},
and the terms are in $E,$
whence $T (a) - L_1 (a) \in E.$
Therefore $\psi (a) = \pi (L_1 (a)) = \rh (A).$
\end{proof}

\begin{cor}\label{P:RelMVb}
Let $A$ be a separable \ca, and identify $A$ with its
image in $\Asym{A}.$
Let $B \subset \Asym{A}$
be a separable subalgebra which contains $A.$
Let $E$ be a nonunital \ca\  with a countable approximate identity.
Then any homomorphism $\ph \colon A \to M (E)/E$ extends to a
homomorphism from $B$ into $M (E)/E.$
\end{cor}

\begin{proof}
The \hm\  $\ph$ induces, in an obvious way,
a \hm\  %
\[
{\overline{\ph}} \colon \Asym{A} \to \Asym{M (E)/E}.
\]
Apply Proposition~\ref{P:RelMV} with $\mu = {\overline{\ph}} |_B.$
\end{proof}

\section{Outer automorphisms of the Calkin algebra}\label{Sec:Pf}

\indent
Most of this section is occupied by the proof of the main theorem.
We start with a useful definition and several related lemmas.
As in the introduction,
we denote by $Q$ the Calkin algebra $L (H) / K (H)$
for a separable infinite dimensional Hilbert space $H.$
We further let $\pi \colon L (H) \to Q$ be the quotient map.

\begin{dfn}\label{D:Split}
Let $A \subset Q$ be a subalgebra,
and let $p \in Q$ be a \pj.
We say that {\emph{$p$ splits $A$}}
if $p$ commutes with every element of $A$ and the
\hm s $x \mapsto p x$ and $x \mapsto (1 - p) x,$
from $A$ to $p Q p$ and to $(1 - p) Q (1 - p),$
are both injective.
We further say that {\emph{$p$ trivially splits $A$}}
if there are a \pj\  $e \in L (H)$ such that $\pi (e) = p$
and a \hm\  $\ph \colon A \to e L (H) e,$
such that $\pi (\ph (x)) = p x$ for all $x \in A.$
\end{dfn}

The additional condition for a trivial splitting
is just that the extension corresponding to $x \mapsto p x$ be trivial.
Note that we do not require the extension corresponding to
$x \mapsto (1 - p) x$ to be trivial.

\begin{lem}\label{L:TrSpl}
Let $A \subset Q$ be a separable unital subalgebra,
and let $l \in \N.$
Then there exist $l$ \mops\  $p_1, p_2, \ldots p_l \in Q,$
each of which trivially splits $A.$
\end{lem}

\begin{proof}
The case $n = 1$ is immediate from the statement that
for any separable unital \ca\  $A,$
any extension by $K (H),$
in the sense of a unital monomorphism $\ta \colon A \to Q,$
is equivalent to its direct sum with a trivial extension
$\ta_0 \colon A \to Q.$
See the discussion at the beginning of Section~4 of~\cite{Ar}.
For the general case,
observe that $\ta \sim \ta \oplus \ta_0$
implies $\ta \sim \ta \oplus \ta_0 \oplus \ta_0,$ etc.
\end{proof}

For commuting \pj s $p$ and $q,$
we write $p \btu q = p + q - 2 p q,$
so that $(1 - 2 p) (1 - 2 q) = 1 - 2 (p \btu q).$

\begin{lem}\label{L:SCnj}
Let $A \subset Q$ be a separable unital subalgebra,
let $p \in Q$ be a \pj\  which splits $A,$
and let $q \in Q$ be a \pj\  which trivially splits
the \ca\  $C^* (A, p)$ generated by $A$ and $p.$
Then there is a norm continuous path of unitaries
$s \mapsto w (s)$ in $Q,$
defined for $s \in [0, 1],$
such that $w (0) = 1,$
$w (1)^* p w (1) = p \btu q,$
$\| w (s_1) - w (s_2) \| \leq \pi | s_1 - s_2 |$
for $s_1, s_2 \in [0, 1],$
and $w (s) x = x w (s)$ for all $s \in [0, 1]$ and $x \in A.$
\end{lem}

\begin{proof}
Write $Q = L (H) / K (H)$ for a separable infinite dimensional
Hilbert space $H.$
Let $D = C^* (A, p) \cong A \oplus A$.
Since $q$ trivially splits $D$ there is a monomorphism
$\ph \colon D \to L (H)$ such that $\pi \circ \ph$ is the map
$x \mapsto q x$ from $D$ to $Q$.
Set $e_0 = \ph (p)$ and $e_1 = \ph (1 - p).$
Then the map
$x \mapsto \ph_0 (x) = \ph (p x)$ is a monomorphism from $A$ to $L (H),$
such that $\ph_0 (A) \subset L (e_0 H) \subset L (H),$
such that $\ph_0 (A) \cap K (H) = \{ 0 \},$
and such that $(\pi \circ \ph_0) (x) = p q x$ for all $x \in A.$
Similarly, $\ph_1 (x) = \ph ((1 - p) x)$
defines a monomorphism $\ph_1 \colon A \to L (e_1 H) \subset L (H)$
such that $(\pi \circ \ph_1) (x) = q (1 - p) x$ for all $x \in A.$
Note that $\ph (x) = \ph_1 (x) + \ph_2 (x)$ for all $x \in A.$
We apply Theorem~5 of~\cite{Ar}.
The \hm s $\ph_0$ and $\ph_1$ are unital
when the codomains are taken to be $L (e_0 H)$ and $L (e_1 H),$
and satisfy Condition~(3) there.
So Condition~(1) there provides, in particular,
an isomorphism $v \colon e_0 H \to e_1 H,$
which we treat as a partial isometry in $L (H),$
such that $v \ph_0 (x) v^* - \ph_1 (x) \in K (H)$ for all $x \in A.$

For $0 \leq s \leq 1,$ define a unitary $z_t \in L (H)$ by
\[
z (s) = (1 - e_0 - e_1)
   + \cos (\pi s) (e_0 + e_1)
   + \sin (\pi s) (v - v^*). 
\]
Set $w (s) = \pi (z (s)).$
One easily checks that $z (s) \ph (x) - \ph (x) z (s) \in K (H)$
for all $x \in A,$
from which one easily concludes that $w (s)$ commutes with every
element of $a.$
{}From
\[
w (1)^* (1 - q) p w (1) = (1 - q) p
\andeqn
w (1)^* q p w (1) = q (1 - p),
\]
we get $w (1)^* p w (1) = p \btu q.$
The other required properties of $s \mapsto w (s)$ are easily verified.
\end{proof}

\begin{thm}\label{T:Main}
Assume the continuum hypothesis.
Then there are $2^{\aleph_1}$ outer automorphisms of $Q.$
\end{thm}



\begin{proof}
Let $(x_{\af})_{\af < \aleph_1}$ be an enumeration of $Q$ by
countable ordinals.
For each countable ordinal $\bt$ and each function
$\ep \colon [0, \bt) \to \{ 0, 1\},$
we construct:
\begin{itemize}
\item
Separable unital subalgebras
$B_{\bt}^{\ep} \subset A_{\bt}^{\ep} \subset Q.$
\item
A \pj\  $p_{\bt}^{\ep} \in A_{\bt}^{\ep}$ and a selfadjoint
unitary $u_{\bt}^{\ep} \in A_{\bt}^{\ep}.$
\item
For every $\gm \in [0, \bt),$
a norm continuous path of unitaries
$s \mapsto w_{\gm, \bt}^{\ep} (s)$ in $A_{\bt}^{\ep},$
defined for $s \in [0, 1].$
\end{itemize}
We will require that these satisfy appropriate conditions,
but, before stating them, we introduce some convenient notation
for functions from ordinals to $\{0, 1 \}.$

If $\ep \colon [0, \bt) \to \{ 0, 1\}$ is a function,
and $\gm < \bt,$
then we write $\ep_{\gm} = \ep |_{[0, \gm)},$
and we frequently even omit the subscript $\gm,$
writing $A_{\gm}^{\ep}$ for $A_{\gm}^{\ep_{\gm}},$
writing $p_{\gm}^{\ep}$ for $p_{\gm}^{\ep_{\gm}},$
and writing $u_{\gm}^{\ep}$ for $u_{\gm}^{\ep_{\gm}}.$
(One may think of
the $A_{\bt}^{\ep}$ as being indexed by
functions $\ep \colon [0, \aleph_1) \to \{ 0, 1 \},$
such that $A_{\bt}^{\ep}$ actually only depends on $\ep |_{[0, \bt)}.$)
Furthermore, if $\ep \colon [0, \bt) \to \{ 0, 1 \}$
and $j \in \{ 0, 1\},$
then we write $\ep \tie j$ for the function
from $[0, \, \bt + 1)$ to $\{ 0, 1\}$
whose restriction to $[0, \bt)$ is $\ep$
and which takes the value~$j$ at $\bt.$

Using this notation, the objects above are required to
satisfy the following conditions for every countable ordinal~$\bt$
and every $\ep \colon [0, \bt) \to \{ 0, 1 \}$:
\begin{enumerate}
\item\label{Ind:Exh}
$x_{\bt} \in A_{\bt}^{\ep}.$
\item\label{Ind:Spl}
$p_{\bt}^{\ep}$ splits $B_{\bt}^{\ep}.$
\item\label{Ind:New}
There exists a unitary $u \in B_{\bt}^{\ep}$
such that $u_{\bt}^{\ep} = u (1 - 2 p_{\bt}^{\ep}).$
If $\bt = \gm + 1,$
then $u = u_{\gm}^{\ep}.$
\item\label{Ind:Inc}
If $\gm < \bt,$ then $A_{\gm}^{\ep} \subset B_{\bt}^{\ep}.$
\item\label{Ind:Aut}
If $\gm < \bt,$ then
$u_{\gm}^{\ep}$ commutes with $u_{\bt}^{\ep}$ and
$( u_{\gm}^{\ep} )^* x u_{\gm}^{\ep}
      = ( u_{\bt}^{\ep} )^* x u_{\bt}^{\ep}$
for all $x \in A_{\gm}^{\ep}.$
\item\label{Ind:Path}
If $\gm < \bt,$
then the path $s \mapsto w (s) = w_{\gm, \bt}^{\ep} (s)$
satisfies
$w (1)^* u_{\gm}^{\ep} w (1) = u_{\bt}^{\ep},$
$w (0) = 1,$
$\| w (s_1) - w (s_2) \| \leq \pi | s_1 - s_2 |$
for $s_1, s_2 \in [0, 1],$
and $w (s) x = x w (s)$ for all $x \in B_{\gm}^{\ep}.$
\item\label{Ind:Dis}
If $\bt = \gm + 1$
and $\ep \colon [0, \gm ) \to \{ 0, 1 \},$
then $A_{\bt}^{\ep \tie 0} = A_{\bt}^{\ep \tie 1}$
and there is $x \in A_{\bt}^{\ep \tie 0}$ such that
$\big( u_{\bt}^{\ep \tie 0} \big)^* x u_{\bt}^{\ep \tie 0}
      \neq \big( u_{\bt}^{\ep \tie 1} \big)^* x u_{\bt}^{\ep \tie 1}.$
\setcounter{TmpEnumi}{\value{enumi}}
\end{enumerate}

We point out that,
for fixed $\ep_0, \ep_1 \colon [0, \aleph_1) \to \{ 0, 1\},$
the subalgebras $A_{\af}^{\ep_0}$ and $A_{\af}^{\ep_1}$ only
start to differ at the first limit ordinal not less than the first
ordinal at which $\ep_0$ and $\ep_1$ disagree.
Unfortunately,
we did not manage to make these subalgebras
fully independent of $\ep.$

Assume the construction has been carried out.
Let $\ep \colon [0, \aleph_1) \to \{ 0, 1\}$ be a function.
As before,
we write $A_{\af}^{\ep}$ for $A_{\af}^{\ep_{\af}},$ etc.
First, Condition~(\ref{Ind:Exh}) implies that
$\bigcup_{\af < \aleph_1} A_{\af}^{\ep} = Q.$
Moreover,
Condition~(\ref{Ind:Aut})
implies that if $x \in A_{\bt}^{\ep},$
then the function $\af \mapsto ( u_{\af}^{\ep} )^* x u_{\af}^{\ep}$
is constant on $[\bt, \aleph_1).$
Therefore
$\ph^{\ep} (x)
 = \lim_{\af \to \aleph_1} ( u_{\af}^{\ep} )^* x u_{\af}^{\ep}$
exists for all $x \in Q.$
Moreover, $\ph^{\ep} \circ \ph^{\ep} = \id_{Q}$ because
$u_{\af}^{\ep}$ is selfadjoint for all $\af.$
So $\ph^{\ep}$ is an automorphism of $Q.$

Now let $\ep_0, \ep_1 \colon [0, \aleph_1) \to \{ 0, 1\},$
with $\ep_0 \neq \ep_1.$
We prove that $\ph^{\ep_0} \neq \ph^{\ep_1}.$
Let $\af$ be the least ordinal such that $\ep_0 (\af) \neq \ep_1 (\af).$
Let $\ep$ be the common restriction
of $\ep_0$ and $\ep_1$ to $[0, \af).$
\Wolog\  $\ep_0 (\af + 1) = 0$
and $\ep_1 (\af + 1) = 1$.
Condition~(\ref{Ind:Dis}) implies that
$A_{\af + 1}^{\ep_0} = A_{\af + 1}^{\ep_1},$
and that there is $x \in A_{\af + 1}^{\ep_0}$ such that
$\big( u_{\af + 1}^{\ep_0} \big)^* x u_{\af + 1}^{\ep_0}
      \neq \big( u_{\af + 1}^{\ep_1} \big)^* x u_{\af + 1}^{\ep_1}.$
It follows from Condition~(\ref{Ind:Aut})
that
\[
\ph^{\ep_0} (x)
 = \big( u_{\af + 1}^{\ep_0} \big)^* x u_{\af + 1}^{\ep_0}
\andeqn
\ph^{\ep_1} (x)
 = \big( u_{\af + 1}^{\ep_1} \big)^* x u_{\af + 1}^{\ep_1}.
\]
Therefore $\ph^{\ep_0} \neq \ph^{\ep_1}.$

Since there are $2^{\aleph_1}$ functions from
$[0, \aleph_1)$ to $\{ 0, 1\},$
it follows that there are at least $2^{\aleph_1}$ automorphisms of $Q.$
However,
there are only $\aleph_1$ unitaries in $Q.$
So there must be $2^{\aleph_1}$ outer automorphisms of $Q.$

We now construct $A_{\af}^{\ep},$ $B_{\af}^{\ep},$
$p_{\af}^{\ep},$ $u_{\af}^{\ep},$
and the $w_{\bt, \gm}^{\ep},$ by transfinite induction on $\af.$

We first consider the case $\af = 0.$
In this case, $[0, \af) = \varnothing,$
and there is a unique function $\ep_0 \colon [0, \af) \to \{ 0, 1\}.$
We let $p_0^{\ep_0} \in Q$ be any nontrivial \pj.
Set $u_0^{\ep_0} = 1 - 2 p_0^{\ep_0}.$
Take $B_0^{\ep_0} = \C,$
and let $A_0^{\ep_0}$ be the subalgebra of $Q$
generated by $1,$ $x_0,$
$p_0^{\ep_0},$ and $u_0^{\ep_0}.$
Since $[0, \af) = \varnothing,$
we do not need any paths $t \mapsto w_{\gm, \bt}^{\ep} (t).$
Condition~(\ref{Ind:Exh}) is satisfied by construction,
Condition~(\ref{Ind:Spl}) is trivially satisfied,
Condition~(\ref{Ind:New}) is satisfied with $u = 1,$
and the remaining conditions are vacuous.

We now consider the induction step.
Assume that $A_{\bt}^{\ep},$ $B_{\bt}^{\ep},$
$p_{\bt}^{\ep},$ and $u_{\bt}^{\ep}$
have been constructed for all $\bt < \af$
and all $\ep \colon [0, \bt) \to \{ 0, 1 \},$
and that the $w_{\bt, \gm}^{\ep}$
have been constructed for all $\gm < \bt < \af$
and all $\ep \colon [0, \bt) \to \{ 0, 1 \}.$

Suppose first that $\af = \bt + 1$ is a successor ordinal.
Let $\ep \colon [0, \bt) \to \{ 0, 1 \}$ be given;
we construct the objects
indexed by $\af$ and by $\ep \tie 0$ and $\ep \tie 1.$
By Lemma~\ref{L:TrSpl} (with $l = 3$),
there are orthogonal \pj s $q_0, q_1 \in Q$
which split $A_{\bt}^{\ep}$ trivially, and such that
$1 - q_0 - q_1 \neq 0.$
By Lemma~\ref{L:SCnj}
there exist norm continuous paths of unitaries
$s \mapsto w_j (s)$ in $Q,$
defined for $j = 0, 1$ and $s \in [0, 1],$
such that $w_j (0) = 1,$
$w_j (1)^* p_{\bt}^{\ep} w_j (1) = p_{\bt}^{\ep} \btu q_j,$
$\| w_j (s_1) - w_j (s_2) \| \leq \pi | s_1 - s_2 |$
for $s_1, s_2 \in [0, 1],$
and $w_j (s) x = x w_j (s)$ for all $x \in B_{\bt}^{\ep}.$
Choose $c \in Q$ such that $c^* c = q_0$ and $c c^* = 1 - q_0 - q_1.$

For $j = 0, 1,$
set $p_{\af}^{\ep \tie j} = q_j,$
set $u_{\af}^{\ep \tie j} = u_{\bt}^{\ep} ( 1 - 2 q_j),$
and set $B_{\af}^{\ep \tie j} = A_{\bt}^{\ep}.$
Let $A_{\af}^{\ep \tie 0} = A_{\af}^{\ep \tie 1}$ be the subalgebra
of $Q$ generated by
$A_{\bt}^{\ep},$ $q_0,$ $q_1,$
$x_{\af},$ $c,$ and $w_j (s)$
for $0 \leq s \leq 1$ and $j = 0, 1.$
Further set
$w_{\bt, \af}^{\ep \tie j} (s) = w_j (s)$
and, for $\gm < \bt,$ set
$w_{\gm, \af}^{\ep \tie j} (s)
   = w_j (1)^* w_{\gm, \bt}^{\ep} (s) w_j (1).$

We check the conditions.
First, $u_{\af}^{\ep \tie j}$ is selfadjoint because
$u_{\bt}^{\ep}$ is selfadjoint and commutes with $q_j.$
Conditions (\ref{Ind:Exh}), (\ref{Ind:Spl}),
and~(\ref{Ind:New})
are satisfied by construction.
Condition~(\ref{Ind:Inc}) is immediate.

We verify Condition~(\ref{Ind:Aut}), for $\gm < \af.$
If $\gm = \bt,$
we observe that $u_{\af}^{\ep \tie j}$
commutes with $u_{\bt}^{\ep \tie j} = u_{\bt}^{\ep}$
because $u_{\bt}^{\ep}$ commutes with $q_j.$
Also, conjugation by both unitaries induces the same map
on $A_{\bt}^{\ep \tie j} = A_{\bt}^{\ep}$
because $q_j$ commutes with every element of $A_{\bt}^{\ep}.$
So suppose that $\gm < \bt.$
Then $u_{\gm}^{\ep \tie j} = u_{\gm}^{\ep}$
commutes with $u_{\bt}^{\ep}$
by the induction hypothesis,
and $u_{\gm}^{\ep}$ commutes with $q_j$
because $u_{\gm}^{\ep} \in A_{\bt}^{\ep},$
so $u_{\gm}^{\ep \tie j}$ commutes with $u_{\af}^{\ep \tie j}.$
Furthermore, conjugation by $u_{\gm}^{\ep \tie j}$ induces the
same map on $A_{\gm}^{\ep \tie j}$
as conjugation by $u_{\bt}^{\ep \tie j}$
by the induction hypothesis,
so also the same map as conjugation by $u_{\af}^{\ep \tie j}$
by the case $\gm = \bt.$

Next, we verify Condition~(\ref{Ind:Path}), for $\gm < \af.$
For $\gm = \bt,$
most of the condition is satisfied by construction.
We must, however, check that
$w_j (1)^* u_{\bt}^{\ep} w_j (1) = u_{\af}^{\ep \tie j}.$
Following~(\ref{Ind:New}) of the induction hypothesis, write
$u_{\bt}^{\ep} = u (1 - 2 p_{\bt}^{\ep})$
for some unitary $u \in B_{\bt}^{\ep}.$
Since $w_j (1)$ commutes with all elements of $B_{\bt}^{\ep},$
we get
\begin{align*}
w_j (1)^* u_{\bt}^{\ep} w_j (1)
 & = w_j (1)^* u (1 - 2 p_{\bt}^{\ep}) w_j (1)
   = u [1 - 2 w_j (1)^* p_{\bt}^{\ep}  w_j (1) ]    \\
 & = u (1 - 2 p_{\bt}^{\ep} \btu q_j)
   = u (1 - 2 p_{\bt}^{\ep}) (1 - 2 q_j)
   = u_{\af}^{\ep \tie j},
\end{align*}
as desired.
Now suppose that $\gm < \bt.$
Set $z (s) = w_{\gm, \af}^{\ep \tie j} (s).$
Then it is immediate from the
definition and the properties of $w_{\gm, \bt}^{\ep} (s)$
assumed in the induction hypothesis that $z (0) = 1$ and
$\| z (s_1) - z (s_2) \| \leq \pi | s_1 - s_2 |$
for $s_1, s_2 \in [0, 1].$
To get $z (s) x = x z (s)$ for $x \in B_{\gm}^{\ep},$
we use in addition the relations
$B_{\gm}^{\ep} \subset B_{\bt}^{\ep}$
and $w_j (s) x = x w_j (s)$ for $x \in B_{\bt}^{\ep}.$
For the remaining part of this condition,
we compute,
using $w_j (1) u_{\gm}^{\ep} w_j (1)^* = u_{\gm}^{\ep}$
and the induction hypothesis at the second step,
\[
z (1)^* u_{\gm}^{\ep} z (1)
 = w_j (1)^* w_{\gm, \bt}^{\ep} (1)^* w_j (1) u_{\gm}^{\ep}
       w_j (1)^* w_{\gm, \bt}^{\ep} (t) w_j (1)
 = w_j (1)^* u_{\bt}^{\ep} w_j (1)
 = u_{\af}^{\ep \tie j}.
\]

It remains to verify Condition~(\ref{Ind:Dis}).
We have $A_{\af}^{\ep \tie 1} = A_{\af}^{\ep \tie 0}$ by construction.
Set $x = u_{\bt}^{\ep} c (u_{\bt}^{\ep})^*,$
which is in $A_{\af}^{\ep \tie 0}.$
Using $c q_1 = q_0 c = q_1 c = 0$ and $c q_0 = c,$
we get
\[
(u_{\af}^{\ep \tie 0})^* x u_{\af}^{\ep \tie 0}
  = (1 - 2 q_0) c (1 - 2 q_0)
  = - c
\]
and
\[
(u_{\af}^{\ep \tie 1})^* x u_{\af}^{\ep \tie 1}
  = (1 - 2 q_1) c (1 - 2 q_1)
  = c.
\]
This completes the proof of the successor ordinal induction step.

Now suppose that $\af$ is a limit ordinal
and $\ep \colon [0, \af) \to \{ 0, 1 \}.$
Choose a strictly increasing sequence of ordinals $(\af_n)_{n \geq 1}$
such that $\lim_{n \to \infty} \af_n = \af.$
Define a \cfn\  $v$ from $[1, \infty)$ to the unitary group of $Q$
by, for $t \in [n, \, n + 1]$ and using~(\ref{Ind:Path}) of the
induction hypothesis,
\[
v (t) = w_{\af_1, \af_2}^{\ep} (1) w_{\af_2, \af_3}^{\ep} (1) \cdots
              w_{\af_{n - 1}, \af_n}^{\ep} (1)
              w_{\af_n, \af_{n + 1}}^{\ep} (t - n).
\]
Further, for each $\bt < \af,$
we define a \cfn\  $c_{\bt}$
from $[0, 1] \times [1, \infty)$ to the unitary group of $Q$
as follows.
Let $m$ be the smallest natural number such that $\bt < \af_m,$
and set
\[
c_{\bt} (s, t)
 = [v (m)^* v (t)]^* w_{\bt, \af_m}^{\ep} (s) [v (m)^* v (t)].
\]

We claim that the objects just defined satisfy the following
properties:
\begin{enumerate}
\setcounter{enumi}{\value{TmpEnumi}}
\item\label{NL:0}
If $n \geq m$ and $t \in [n, \, n + 1],$
then
\[
v (m)^* v (t)
  = w_{\af_m, \af_{m + 1}}^{\ep} (1) \cdots
           w_{\af_{n - 1}, \af_n}^{\ep} (1)
           w_{\af_n, \af_{n + 1}}^{\ep} (t - n).
\]
\item\label{NL:N}
For $n \in \N,$ we have
$v (n)^* u_{\af_1}^{\ep} v (n) = u_{\af_n}^{\ep}.$
\item\label{NL:1}
For $t \geq n + 1$ and $x \in A_{\af_n}^{\ep},$ we have
\[
\big[ v (t)^* u_{\af_1}^{\ep} v (t) \big]^*
    x \big[ v (t)^* u_{\af_1}^{\ep} v (t) \big]
 = ( u_{\af_n}^{\ep} )^* x u_{\af_n}^{\ep}.
\]
\item\label{NL:2}
If $\bt < \af_m$ and $t \geq m,$ then
$c_{\bt} (1, t)^* u_{\bt}^{\ep} c_{\bt} (1, t)
   = v (t)^* u_{\af_1}^{\ep} v (t).$
\item\label{NL:3}
For $\bt < \af$ and $t \in [0, \infty),$
we have $c_{\bt} (0, t) = 1.$
\item\label{NL:3b}
For $\bt < \af,$ $t \in [0, \infty),$ and $s_1, s_2 \in [0, 1],$
we have
$\| c_{\bt} (s_1, t) - c_{\bt} (s_2, t) \| \leq \pi | s_1 - s_2 |.$
\item\label{NL:5}
For $\gm < \bt < \af_m,$ $t \in [m, \infty),$ $s \in [0, 1],$
and $x \in B_{\gm}^{\ep},$
the element $c_{\bt} (s, t)$ commutes with $x.$
\setcounter{TmpEnumi}{\value{enumi}}
\end{enumerate}

Property~(\ref{NL:0}) is immediate from the definition,
and Property~(\ref{NL:N}) follows by induction
from~(\ref{Ind:Path}) of the induction hypothesis.

We prove~(\ref{NL:1}).
Let $x \in A_{\af_n}^{\ep}$ and $t \in [n + k, \, n + k + 1]$
with $k \geq 1.$
Set $c = w_{\af_{n + k}, \, \af_{n + k + 1}}^{\ep} (t - n - k).$
So $v (t)^* u_{\af_1}^{\ep} v (t) = c^* u_{\af_{n + k}}^{\ep} c$
by~(\ref{NL:N}).
Using~(\ref{Ind:Aut}) of the induction hypothesis, we get
$(u_{\af_{n + k}}^{\ep} )^* x u_{\af_{n + k}}^{\ep}
  = (u_{\af_{n}}^{\ep} )^* x u_{\af_{n}}^{\ep},$
and this element is in $A_{\af_{n}}^{\ep}$
because $u_{\af_{n}}^{\ep} \in A_{\af_{n}}^{\ep}.$
It now follows from
(\ref{Ind:Path}) and~(\ref{Ind:Inc}) of the induction hypothesis
that $c$ commutes with $x$ and with
$(u_{\af_{n + k}}^{\ep} )^* x u_{\af_{n + k}}^{\ep}.$
Therefore 
\begin{align*}
\big[ v (t)^* u_{\af_1}^{\ep} v (t) \big]^*
    x \big[ v (t)^* u_{\af_1}^{\ep} v (t) \big]
& = \big[ c^* u_{\af_{n + k}}^{\ep} c \big]^*
    x \big[ c^* u_{\af_{n + k}}^{\ep} c \big]  \\
& = c^* (u_{\af_{n + k}}^{\ep})^* x u_{\af_{n + k}}^{\ep} c
  = (u_{\af_{n}}^{\ep} )^* x u_{\af_{n}}^{\ep},
\end{align*}
as desired.

To prove~(\ref{NL:2}), let $t \in [n, \, n + 1]$ with $n \geq m.$
Property~(\ref{NL:0}) (which we already proved)
and Conditions (\ref{Ind:Path}) and~(\ref{Ind:Inc})
of the induction hypothesis
imply that $v (m)^* v (t)$ commutes with $u_{\bt}^{\ep}.$
Using this at the second step,
(\ref{Ind:Path}) of the induction hypothesis at the third step,
and~(\ref{NL:N}) at the fourth step, we get
\begin{align*}
&
c_{\bt} (1, t)^* u_{\bt}^{\ep} c_{\bt} (1, t)
        \\
& \mbox{} \hspace*{2em}
  = [v (m)^* v (t)]^* w_{\bt, \af_m}^{\ep} (1)^* [v (m)^* v (t)]
      u_{\bt}^{\ep}
    [v (m)^* v (t)]^* w_{\bt, \af_m}^{\ep} (1) [v (m)^* v (t)]
        \\
& \mbox{} \hspace*{2em}
  = [v (m)^* v (t)]^* w_{\bt, \af_m}^{\ep} (1)^*
      u_{\bt}^{\ep} w_{\bt, \af_m}^{\ep} (1) [v (m)^* v (t)]
        \\
& \mbox{} \hspace*{2em}
  = [v (m)^* v (t)]^* u_{\af_m}^{\ep} [v (m)^* v (t)]
  = v (t)^* u_{\af_1}^{\ep} v (t),
\end{align*}
as desired.

Properties (\ref{NL:3}) and~(\ref{NL:3b})  
follow from the definitions
and~(\ref{Ind:Path}) of the induction hypothesis.
For Property~(\ref{NL:5}), one uses in addition
Property~(\ref{NL:0}).
This completes the proof of (\ref{NL:0}) through~(\ref{NL:5}).

Now set $M = {\overline{\bigcup_{n = 0}^{\infty} A_{\af_n}^{\ep}}},$
which of course is equal to
${\overline{\bigcup_{\bt < \af} A_{\bt}^{\ep}}}.$
We have $v (t) \in M$ for all $t$
and $c_{\bt} (s, t) \in M$ for all $s,$ $t,$ and $\bt.$
Identify $M$ with its image in $\Asym{M},$
let ${\overline{v}}$ be the image of $v$ in this algebra,
and let ${\overline{c}}_{\bt} (s)$ be the image of
the function $t \mapsto c_{\bt} (s, t)$ in this algebra.
Property~(\ref{NL:2}) implies that
${\overline{c}}_{\bt} (1)^* u_{\bt}^{\ep} {\overline{c}}_{\bt} (1)
   = {\overline{v}}^* u_{\af_1}^{\ep} {\overline{v}}$
for all $\bt < \af.$

Let $N \subset \Asym{M}$ be the \ca\  generated
by $M$ and ${\overline{v}}.$
Note that ${\overline{c}}_{\bt} (s) \in N$ for all $s \in [0, 1]$
and all $\bt < \af.$
Apply Corollary~\ref{P:RelMVb} with $M$ in place of $A$
and $N$ in place of $B,$
getting a \hm\  $\ph \colon N \to Q$
such that $\ph (x) = x$ for all $x \in M.$
Set $z = \ph ({\overline{v}})$
and set
$y_{\bt} (s) = \ph ({\overline{c}}_{\bt} (s))
 = [v (m)^* z]^* w_{\bt, \af_m}^{\ep} (s) [v (m)^* z]$
for $\bt < \af.$
Further set $u = z^* u_{\af_1}^{\ep} z.$
We claim that these elements satisfy the following properties:
\begin{enumerate}
\setcounter{enumi}{\value{TmpEnumi}}
\item\label{Prp:SA}
$u$ is selfadjoint.
\item\label{Prp:Cmm}
$u$ commutes with $u_{\bt}^{\ep}$ for every $\bt < \af.$
\item\label{Prp:Aut}
For every $\bt < \af$ and every $x \in A_{\bt}^{\ep},$ we have
$u x u^* = (u_{\bt}^{\ep})^* x u_{\bt}^{\ep}.$
\item\label{Prp:One}
$y_{\bt} (0) = 1$ for all $\bt < \af.$
\item\label{Prp:Cnj}
$y_{\bt} (1)^* u_{\bt}^{\ep} y_{\bt} (1) = u$
for all $\bt < \af.$
\item\label{Prp:Lip}
$\| y_{\bt} (s_1) - y_{\bt} (s_2) \| \leq \pi | s_1 - s_2 |$
for all $s_1, s_2 \in [0, 1]$ and $\bt < \af.$
\item\label{Prp:yCm}
$y_{\bt} (s) x = x y_{\bt} (s)$ for all $\bt < \af,$
all $x \in B_{\bt}^{\ep},$ and all $s \in [0, 1].$
\end{enumerate}

Property~(\ref{Prp:SA}) is immediate.
Properties (\ref{Prp:One}), (\ref{Prp:Cnj}),
(\ref{Prp:Lip}), and~(\ref{Prp:yCm})
follow from, in order,
Properties (\ref{NL:3}), (\ref{NL:2}),
(\ref{NL:3b}), and~(\ref{NL:5}).
We prove Properties (\ref{Prp:Cmm}) and~(\ref{Prp:Aut}).
Property~(\ref{NL:1}) implies that
for every $n$ and for every $x \in A_{\af_n}^{\ep},$
we have
$(u_{\af_n}^{\ep})^* x u_{\af_n}^{\ep} = u^* x u.$
For $\gm < \af,$ choose $n$ with $\af_n > \gm$
and use Condition~(\ref{Ind:Aut}) of the induction hypothesis
for~$\af_n$ to get
$(u_{\gm}^{\ep})^* x u_{\gm}^{\ep} = u^* x u$
for all $x \in A_{\gm}^{\ep}.$
Putting $\gm = \bt$ gives Property~(\ref{Prp:Aut}).
For Property~(\ref{Prp:Cmm}),
put $x = u_{\gm}^{\ep}$ and $\gm = \bt + 1,$
and use~(\ref{Ind:Aut}) of the induction hypothesis for~$\gm.$
This completes the proof
of Properties (\ref{Prp:SA}) through~(\ref{Prp:yCm}).

Let $B_{\af}^{\ep}$ be the subalgebra of $Q$
generated by $M$ and $z.$
Use Lemma~\ref{L:TrSpl} to
choose a \pj\  $p_{\af}^{\ep}$ which splits $B_{\af}^{\ep}$
trivially, and
set $u_{\af}^{\ep} = u (1 - 2 p_{\af}^{\ep}).$
For each $\bt < \af,$
the element $u_{\bt}^{\ep}$ commutes with $u$
by Property~(\ref{Prp:Cmm}),
and both are selfadjoint,
so there exists a \pj\  $q_{\bt}$ such that
$u (u_{\bt}^{\ep})^* = (u_{\bt}^{\ep})^* u = 1 - 2 q_{\bt}.$
Then $q_{\bt}$ commutes with every element of $A_{\bt}^{\ep}$
by Property~(\ref{Prp:Aut}).
We claim that $q_{\bt}$ splits $A_{\bt}^{\ep}.$
To see this, recall that by~(\ref{Ind:New}) of the induction hypothesis
we have $u_{\bt + 1}^{\ep} = u_{\bt}^{\ep} (1 - 2 p_{\bt + 1}^{\ep}).$
Using Property~(\ref{Prp:Cnj}) at the second step
and Property~(\ref{Prp:yCm}) at the fourth step, we then get
\begin{align*}
u_{\bt}^{\ep} \big( 1 - 2 q_{\bt} \big)
& = u = y_{\bt+1} (1)^* u_{\bt+1}^{\ep} y_{\bt+1} (1) \\
& = y_{\bt+1} (1)^* u_{\bt}^{\ep}
        \big( 1 - 2 p_{\bt + 1}^{\ep} \big) y_{\bt+1} (1) \\
& = u_{\bt}^{\ep} y_{\bt+1} (1)^*
        \big( 1 - 2 p_{\bt + 1}^{\ep} \big) y_{\bt+1} (1) \\
& = u_{\bt}^{\ep}
    \big[ 1 - 2 y_{\bt+1} (1)^* p_{\bt+1}^{\ep} y_{\bt+1} (1) \big],
\end{align*}
which implies that
$q_{\bt} = y_{\bt + 1} (1)^* p_{\bt + 1}^{\ep} y_{\bt + 1} (1).$
Now $p_{\bt + 1}^{\ep}$ splits $A_{\bt}^{\ep}$ by~(\ref{Ind:Spl})
and~(\ref{Ind:Inc}) of the induction hypothesis.
For any nonzero $x \in A_{\bt}^{\ep}$ we therefore have
(using Property~(\ref{Prp:yCm}))
\[
q_{\bt} x = y_{\bt + 1} (1)^* p_{\bt + 1}^{\ep} y_{\bt + 1} (1) x
= y_{\bt + 1} (1)^* p_{\bt + 1}^{\ep} x y_{\bt + 1} (1) \neq 0,
\]
and $(1 - q_{\bt}) x \neq 0$ by an analogous computation.
So $q_{\bt}$ splits $A_{\bt}^{\ep}.$

Lemma~\ref{L:SCnj}
therefore provides, in particular (ignoring the path)
a unitary $w_{\bt}$ (which would be $w_{\bt} (1)$)
such that $w_{\bt}^* q_{\bt} w_{\bt} = q_{\bt} \btu p_{\af}^{\ep}$
and $w_{\bt} x = x w_{\bt}$ for all $x \in A_{\bt}^{\ep}.$
Writing
$u = u_{\bt}^{\ep} (u_{\bt}^{\ep})^* u
             = u_{\bt}^{\ep} (1 - 2 q_{\bt}),$
and using $u_{\bt}^{\ep} \in A_{\bt}^{\ep},$
we get $w_{\bt}^* u w_{\bt} = u_{\af}^{\ep}.$
Set $w_{\bt, \af} (s) = w_{\bt}^* y_{\bt} (s) w_{\bt}.$
Then let $A_{\af}^{\ep}$ be the subalgebra of $Q$
generated by $B_{\af}^{\ep},$ $x_{\af},$ and $p_{\af}^{\ep}.$

We verify that the required properties hold.
Conditions (\ref{Ind:Exh}), (\ref{Ind:Spl}),
(\ref{Ind:New}), and~(\ref{Ind:Inc}) are clearly satisfied,
while Condition~(\ref{Ind:Dis}) is vacuous.
Condition~(\ref{Ind:Aut}) follows from
Properties (\ref{Prp:Cmm}) and~(\ref{Prp:Aut}),
together with the fact that $p_{\af}^{\ep}$ commutes
with $u$ and with every element
of $\bigcup_{\bt < \af} A_{\bt}^{\ep}.$

It remains only to verify Condition~(\ref{Ind:Path}).
We replace $y_{\bt} (s)$ with $w_{\bt}^* y_{\bt} (s) w_{\bt}$
in each of
Properties (\ref{Prp:One}), (\ref{Prp:Cnj}),
(\ref{Prp:Lip}), and~(\ref{Prp:yCm}).
Properties (\ref{Prp:One}) and (\ref{Prp:Lip})
immediately become the corresponding parts of
Condition~(\ref{Ind:Path}).
Properties (\ref{Prp:Cnj}) and (\ref{Prp:yCm})
become the corresponding parts of Condition~(\ref{Ind:Path})
after using the fact that $w_{\bt}$ commutes with $u_{\bt}^{\ep}$
and every element of $B_{\bt}^{\ep}.$
This completes the limit ordinal induction step,
and hence the proof.
\end{proof}

\begin{rmk}\label{R:AltPf}
If one is only concerned with producing a single outer automorphism
of $Q,$
the following more direct approach can be used.
Let $(\ps_{\af})_{\af < \aleph_1}$ be an enumeration
by countable ordinals of
the inner automorphisms of $Q.$
At each successor ordinal step,
instead of constructing two \pj s
$p_{\af}^{\ep \tie 0}$ and $p_{\af}^{\ep \tie 1}$
such that the corresponding unitaries
$u_{\af}^{\ep \tie 0}$ and $u_{\af}^{\ep \tie 1}$
give different inner automorphisms of
$A_{\af}^{\ep \tie 0} = A_{\af}^{\ep \tie 1},$
choose a single \pj\  $p_{\af},$
but use the existence of more than one choice to
ensure that conjugation by the corresponding unitary $u_{\af}$
disagrees with $\ps_{\af}$ on $A_{\af}.$
This yields an automorphism $\ph$ of $Q$
which is not equal to any $\ps_{\af},$
and hence not inner.

In fact, this argument shows that there are strictly more
than $\aleph_1$ automorphisms of $Q.$
However, from the set-theoretic point of view,
it is a slightly stronger conclusion
that there are $2^{\aleph_1}$ automorphisms of $Q.$
\end{rmk}

\end{document}